# A FAST SOLVER FOR PENTADIAGONAL TOEPLITZ SYSTEMS

TECHNICAL REPORT

**Shahin Hasanbeigi**
Department of Applied Mathematics
Tarbiat Modares University
Tehran, Iran
shahin77hb@gmail.com

April 2024

**ABSTRACT**

The objective of this work is to present a novel approach for the solution of Pentadiagonal Toeplitz systems of equations that is both faster and more effective than existing classical direct methods. The distinctive structure of Pentadiagonal Toeplitz matrices can be leveraged to devise an algorithm for solving upper triangle systems, rather than the original system. This approach is considerably more straightforward and expeditious than classical methods such as LU and Gauss Eliminations. A comparison with the LU and PLU methods demonstrates the efficacy of our novel algorithm. Furthermore, numerical tests substantiate this efficacy.

***Keywords*** pentadiagonal Toeplitz Systems · Fast Solver · Block LU Factorization · Toeplitz Matrices

## 1 Introduction

Consider the following nonsingular linear system of equation

$$Ax = b \tag{1}$$

where A is an $n \times n$ pentadiagonal Toeplitz matrix of the form

$$A = \begin{bmatrix} \alpha & \beta & \gamma & & & & & & \\ \lambda & \alpha & \beta & \gamma & & & & & \\ \sigma & \lambda & \alpha & \beta & \gamma & & & & \\ & \ddots & \ddots & \ddots & \ddots & \ddots & & & \\ & & \ddots & \ddots & \ddots & \ddots & \ddots & & \\ & & & \ddots & \ddots & \ddots & \ddots & \ddots & \\ & & & & \sigma & \lambda & \alpha & \beta & \gamma \\ & & & & & \sigma & \lambda & \alpha & \beta \\ & & & & & & \sigma & \lambda & \alpha \end{bmatrix} \tag{2}$$

Toeplitz matrices represent a distinctive class of structured matrices that have been employed in a multitude of applications within the domains of applied mathematics, engineering, and scientific computing. The unique structure of Toeplitz matrices has been the subject of numerous research papers and works in the field of linear system solvers, see [1, 2].

In general, there are two principal approaches to the solution of Toeplitz linear systems. Two principal approaches are available for the solution of Toeplitz linear systems: direct methods and iterative methods. Iterative methods are primarily comprised of classical splitting iteration methods and Krylov subspace iteration methods. Classical splitting iteration methods for Toeplitz systems necessitate the implementation of efficient splittings, which are contingent upon the structural and intrinsic characteristics of the coefficient matrices. Illustrative examples of such splittings include Gauss-Seidel and SOR splittings for H-matrices and Hermitian positive definite matrices, Circulant and Skew-Circulant splittings for positive definite matrices, and other forms of splitting iteration methods tailored to specific structural and intrinsic characteristics of the coefficient matrices. Krylov subspace iteration methods include the CG method for Hermitian positive definite matrices and the GMRES method for non-Hermitian positive definite matrices, as well as their variants, as cited in Chan et al. [1, 5, 3] In contrast, direct methods, as described in [2], are often applicable to small and moderate-sized problems but are often too expensive to be practical for large sparse problems. Moreover, direct methods may be susceptible to numerical instability and the loss of solution accuracy, [5].

Pentadiagonal Toeplitz matrices have been identified as playing a role in a multitude of mathematical, scientific, and engineering investigations. As a case in point, pentadiagonal Toeplitz matrices emerge in the context of second and fourth-order differential equations, with different boundary conditions, [6]. The special structure of pentadiagonal Toeplitz matrices can be exploited to develop a rapid algorithm for solving systems of equations with pentadiagonal Toeplitz matrices as their coefficient matrix.

This paper is organized as follows: In the following section, we develop some fast algorithms for solving (1) with the coefficient matrix A being a pentadiagonal Toeplitz matrix, analogous to the matrix in (2). The efficacy of our algorithms is demonstrated through numerical tests in Section 3, after which a brief conclusion is drawn in Section 4.

## 2 Fast Algorithm

Let

$$J = \begin{bmatrix} 0 & 0 & 1 & 0 & 0 & \cdots & 0 \\ 0 & 0 & 0 & 1 & 0 & \cdots & 0 \\ \vdots & \vdots & \vdots & \vdots & \ddots & \ddots & \vdots \\ 0 & 0 & 0 & 0 & \cdots & \cdots & 1 \\ 1 & 0 & 0 & 0 & \cdots & \cdots & 0 \\ 0 & 1 & 0 & 0 & \cdots & \cdots & 0 \end{bmatrix}$$

It is easy to verify that $\hat{A} = JA$ has the following $3 \times 3$ block form

$$\hat{A} = \left[\begin{array}{ccccccc|c|c} \sigma & \lambda & \alpha & \beta & \gamma & & & & \\ & \ddots & \ddots & \ddots & \ddots & \ddots & & & \\ & & \ddots & \ddots & \ddots & \ddots & \ddots & & \\ & & & \ddots & \ddots & \ddots & \ddots & \gamma & \\ & & & & \sigma & \lambda & \alpha & \beta & \gamma \\ & & & & & \sigma & \lambda & \alpha & \beta \\ & & & & & & \sigma & \lambda & \alpha \\ \hline \alpha & \beta & \gamma & & & & & 0 & 0 \\ \hline \lambda & \alpha & \beta & \gamma & & & & 0 & 0 \end{array}\right] \equiv \left[\begin{array}{c|c|c} A_{11} & p & r \\ \hline w^T & 0 & 0 \\ \hline s^T & 0 & 0 \end{array}\right] \tag{3}$$

and has the $3 \times 3$-LU factorization

$$\hat{A} = \begin{bmatrix} I & 0 & 0 \\ w^T A_{11}^{-1} & 1 & 0 \\ s^T A_{11}^{-1} & \frac{s^T A_{11}^{-1} p}{w^T A_{11}^{-1}} & 1 \end{bmatrix} \begin{bmatrix} A_{11} & p & r \\ 0 & -w^T A_{11}^{-1} p & -w^T A_{11}^{-1} r \\ 0 & 0 & -s^T A_{11}^{-1} r - \frac{s^T A_{11}^{-1} p}{w^T A_{11}^{-1} p} w^T A_{11}^{-1} r \end{bmatrix} \tag{4}$$

By multiplying $J$ to both sides of equation (1), we have

$$\hat{A}x = \hat{b} \tag{5}$$

where $\hat{b} = (b_3, b_4, ..., b_n, b_1, b_2)^T$.

Next, we partitioning $x$ and $\hat{b}$ into the following forms

$$x = \left[\begin{array}{c} x_1 \\ \hline x_{n-1} \\ \hline x_n \end{array}\right] \quad , \quad \hat{b} = \left[\begin{array}{c} b_3 \\ \hline b_1 \\ \hline b_2 \end{array}\right] \tag{6}$$

where $x_1$ and $b_3$ are $(n-2) \times 1$ vectors, and $x_{n-1}, x_n, b_1$ and $b_2$ are scalers.

By using the block LU-factorization (4), we have

$$\begin{cases} A_{11}x_1 + x_{n-1}p + x_x r & = b_3 \\ w^T A_{11}^{-1} b_3 - w^T A_{11}^{-1} p x_{n-1} - w^T A_{11}^{-1} r x_n & = b_1 \\ s^T A_{11}^{-1} b_3 - s^T A_{11}^{-1} p x_{n-1} - s^T A_{11}^{-1} x_n & = b_2 \end{cases} \tag{7}$$

For solving (7), we first solve

$$\begin{cases} A_{11}u & = b_3 \\ A_{11}v & = p \\ A_{11}z & = r \end{cases} \tag{8}$$

by the following algorithm

**Algorithm 1** Backward Substitution for Solving $A_{11}y = c$

**Require:** $\sigma, \lambda, \alpha, \beta, \gamma$ and **c**

1: Let $\sigma_1 = \frac{1}{\sigma}, \alpha_1 = \frac{\alpha}{\sigma}, \beta_1 = \frac{\beta}{\sigma}, \gamma_1 = \frac{\gamma}{\sigma}, \lambda_1 = \frac{\lambda}{\sigma}$

2: Compute :

$$\begin{aligned} y_n &= \sigma_1 c_n \\ y_{n-1} &= \sigma_1 c_{n-1} - \lambda_1 y_n \\ y_{n-2} &= \sigma_1 c_{n-2} - \lambda_1 y_{n-1} - \alpha_1 y_n \\ y_{n-3} &= \sigma_1 c_{n-3} - \lambda_1 y_{n-2} - \alpha_1 y_{n-1} - \beta_1 y_n \end{aligned}$$

3: **for** $k = n-4, ..., 1$ **do**

Compute : $y_k = \sigma_1 c_k - \lambda_1 y_{k+1} - \alpha_1 y_{k+2} - \beta_1 y_{k+3} - \gamma_1 y_{k+4}$

4: **end for**

5: **return** $y = [y_1, ..., y_n]^T$

Using algorithm (1), we can solve (8). After we get the solutions $u$ , $v$ and $z$ from (8), the solution of (1) turns into the form

$$\begin{cases} x_1 + x_{n-1}v + x_n r & = u \\ w^T u - w^T v x_{n-1} - w^T z x_n & = b_1 \\ s^T u - s^T v x_{n-1} - s^T z x_n & = b_2 \end{cases} \tag{9}$$

Second and third equations are free from $x_1$ and form a simple $2 \times 2$ system of equations for $x_{n-1}$ and $x_n$ as follow

$$\begin{cases} w^T v x_{n-1} + w^T z x_n & = w^T u - b_1 \\ s^T v x_{n-1} + s^T z x_n & = s^T u - b_2 \end{cases} \tag{10}$$

or its equivalent matrix form

$$\begin{bmatrix} w^T u & w^T z \\ s^T u & s^T z \end{bmatrix} \begin{bmatrix} x_{n-1} \\ x_n \end{bmatrix} = \begin{bmatrix} w^T u - b_1 \\ s^T u - b_2 \end{bmatrix} \tag{11}$$

Since all elements in above equations are scalers, solving (10) or (11) is very easy and not expensive. It's easy to see that

$$\begin{cases} x_{n-1} & = \frac{w^T u - b_1}{w^T v} - \frac{w^T z}{w^T z} \frac{s^T u - b_2 - \frac{s^T v}{w^T v}(w^T u - b_1)}{s^T z - \frac{s^T v}{w^T v} w^T z} \\ \\ x_n & = \frac{s^T u - b_2 - \frac{s^T v}{w^T v}(w^T u - b_1)}{s^T z - \frac{s^T v}{w^T v} w^T z} \end{cases} \tag{12}$$

then, solution of (1) is

$$\begin{cases} x_1 & = u - x_{n-1}v - x_n r \\ x_{n-1} & = \frac{w^T u - b_1}{w^T v} - \frac{w^T z}{w^T z}\frac{s^T u - b_2 - \frac{s^T v}{w^T v}(w^T u - b_1)}{s^T z - \frac{s^T v}{w^T v} w^T z} \\ x_n & = \frac{s^T u - b_2 - \frac{s^T v}{w^T v}(w^T u - b_1)}{s^T z - \frac{s^T v}{w^T v} w^T z} \end{cases} \tag{13}$$

By above analysis we can write our algorithm as follow

**Algorithm 2** Fast Algorithm for Pentadiagonal Toeplitz systems

**Require:** Matrix $A_{11}$ , vectors $p, r, w^T$ and $s^T$ from $\hat{A}$ ,
and vector $b_3$ and scalers $b_2$ and $b_1$ from $\hat{b}$.
1: Solve (8) with algorithm (1) for $u$,$v$ and $z$ .
2: Compute $x_n$ , $x_{n-1}$ and $x_1$ according to (13).
3: **return** $x = [x_1^T \;\; x_{n-1} \;\; x_n]^T$.

In algorithm (2), when compting $x_{n-1}$ , $x_n$ and $x_1$, $u$ , $v$ and $r$ are used twice, $w^T u$ , $w^T v$ , $w^T z$ , $s^T u$ , $s^T v$ and $s^T z$ are used once, so influence of error propagation is very small, therefore computed $x_1$ , $x_{n-1}$ and $x_n$ are reliable. Thus we can conclude that algorithm (1) for solving Pentadiagonal Toeplitz linear systems is numerically stable and computed answer is reliable.

## 3 Numerical Experiments

In this section we use some examples to show effectiveness of algorithm (1) in compair to other existing algorithms.

All the numerical tests were done on an ASUS laptop PC with AMD A12 CPU, 8Gb RAM and by Matlab R2016(b) with a machine precision of $10^{-16}$. For convenience, throughout our numerical experiments, we denote by Relative error $= ||b - Ax||/||b||$, the relative residual error, computing time (in seconds), our new algorithm, LU factorization method with pivoting and LU factorization method. In all tables, the Time is the average value of computing times required by performing the corresponding algorithm 10 times, the right hand side vector b is taken to be $Ax^*$ with $x^* = rand(n, 1)$.

### Examples :

Some artificial examples were used in this experiment. The values of $\alpha, \beta, \gamma, \lambda$ and $\sigma$ corresponding to each examples are presented in Table(1). Results for solving corresponding tests from Table(1) are shown in Tables (2,3,4).

Table 1: Values for $\alpha, \beta, \gamma, \lambda$ and $\sigma$

| | $\sigma$ | $\lambda$ | $\alpha$ | $\beta$ | $\gamma$ |
|---|---|---|---|---|---|
| **Test 1** | 5 | 2 | 4 | 1 | 3 |
| **Test 2** | 1 | 0.2 | 0.1 | 0.2 | 0.5 |
| **Test 3** | 28 | 19 | 17 | 21 | 25 |

Table 2: Relative errorr and time for Test 1

| Test 1 | | n=$2^7$ | n=$2^8$ | n=$2^9$ |
|---|---|---|---|---|
| Relative error | New | 9.7873 e -17 | 1.0744 e -16 | 5.1306 e -14 |
| | PLU | 0.0629 | 0.0443 | 145.7111 |
| Time | New | 1.7919 e -4 | 1.6492 e -4 | 1.4548 e -4 |
| | PLU | 0.0024 | 0.0058 | 0.0109 |

Table 3: Relative errorr and time for Test 2

| Test 2 | | n=$2^7$ | n=$2^8$ | n=$2^9$ |
|---|---|---|---|---|
| Relative error | New | 1.1445 e -16 | 1.1600 e -16 | 0.0167 |
| | PLU | 0.1083 | 0.0769 | 246.7 e 13 |
| Time | New | 1.3660 e -4 | 1.5174 e -4 | 1.9009 e -4 |
| | PLU | 0.0025 | 0.0049 | 0.0276 |

Table 4: Relative errorr and time for Test 3

| Test 3 | | n=$2^7$ | n=$2^8$ | n=$2^9$ |
|---|---|---|---|---|
| Relative error | New | 1.1872 e -16 | 1.0819 e -16 | 1.1157 e -15 |
| | PLU | 0.0656 | 0.0432 | 0.0326 |
| Time | New | 1.4540 e -4 | 1.2301 e -4 | 2.0523 e -4 |
| | PLU | 0.0025 | 0.0049 | 0.0109 |

## 4 Conclusion

In this work, we introduce a fast algorithm for solving pentadiagonal Toeplitz linear systems. This novel approach represents a direct extension of the methodology presented in [7] for tridiagonal Toeplitz systems, exhibiting a high degree of similarity in its computational procedures. We initially exploited the special structure of the Pentadiagonal Toeplitz matrix and reduced the primary computations to the solution of three upper triangular systems in (8), which proved to be considerably more efficient than the original system.